# The Koopmanization of controlled nonlinear Itô stochastic differential systems and its comparison with the Carleman embedding: new results

Amruta Lambe[a], Shambhu Nath Sharma[b]

Department of Electrical Engineering,
Sardar Vallabhbhai National Institute of Technology Surat 395007, Gujarat, India

[a]amrutalambe121@gmail.com, [b]sns@eed.svnit.ac.in

**Abstract:** The Koopmanization embeds the bilinearization via the action of the infinitesimal stochastic Koopman operator on the observables associated with the controlled nonlinear Itô stochastic differential system without explicit linearizations. The stochastic evolutions of controlled Markov processes assume the structure of controlled nonlinear Itô stochastic differential equations. This paper sketches a Koopman operator framework for the filtering of the controlled nonlinear Itô stochastic differential system. The major ingredients of this paper are the construction of the eigenfunctions, action of the infinitesimal stochastic Koopman operator, multi-dimensional Itô differential rule and filtering concerning the controlled nonlinear Itô stochastic differential system. In this paper, we illustrate the 'filtering in the Koopman setting' for a polynomial system and compare with the filtering in the Carleman setting.

## 1. Introduction

The paper analyses the controlled nonlinear Itô stochastic differential equation (SDE) of the form

$$dx_t = f(x_t, u_t)dt + \sum_{1 \leq \gamma \leq r} g_\gamma(x_t, u_t)dB_\gamma,$$

given the observation equation $dy_t = h(x_t)dt + r_1 d\eta_t$ in the Koopman operator framework, where the process $x_t$ is a controlled Markov process. Then, we design filtering algorithms of the concerning system, where the notations have standard meanings. Note that $B_t$ and $\eta_t$ are the independent standard Brownian motion processes.

"The $r-$dimensional random process $B_t$ acting on 'the nonlinear dynamic system with a $d-$dimensional control parameter $u_t$' that results in the controlled Markov process trajectory" is expressible using the above setup of stochastic differential equation. In stochastic control, the explicitly computable control of the original nonlinear stochastic differential system can be achieved via the linearization of the nonlinearity and recasting the associated stochastic differential equation with linear in the control parameter (Kushner, 1967). The nonlinear filtering (Pugachev & Sinitsyn, 1987; Sharma, 2009) is a potential problem for the controlled nonlinear stochastic system. That can be designed by exploiting an equivalent 'no control' setting of the given original system. Other methods are available in literature to achieve 'nonlinear filtering' of the controlled stochastic differential system, which is input to the feedback for the construction of the control parameter using the certainty equivalence principle (James, 1994). The Carleman linearization (Carleman, 1932; Kowalski & Steeb, 1991; Belabbas & Chen, 2023) has the property to allow the application of the generalized linear system theory to controlled nonlinear systems in the deterministic and stochastic frameworks. Besides the methods and techniques available in stochastic control literature, Germani *et al.* (2007) developed an alternative and an appealing unified theory of nonlinear filtering of a controlled nonlinear Itô SDE by using the applicability of the Carleman linearization. The analytic properties of the drift and process noise coefficient allow the applicability of the Carleman linearization. The Carleman linearization introduces the bilinearization into the original controlled nonlinear differential system with respect to the associated linear forcing term. The bilinearization simplifies the setup of the original controlled nonlinear system and culminates into the bilinear form (Bhatt & Sharma, 2020) as well as allows the application of the generalized linear system-theoretic framework. For the bilinearization of the original controlled nonlinear system, the Koopman operator framework has found its potential usefulness as well (Mauroy, Mezić & Sususki, 2020; Brunton, Budišić, Kaiser & Kutz, 2022). In Mauroy *et al.* (2020, p.249), the Koopman operator framework was shown as an alternative to the Carleman linearization with respect to the controlled nonlinear systems. In the sense of achieving 'greater refinements' in the linearization of nonlinear systems, the Carleman linearization and the Koopmanization can be termed as the 'super-linearization' techniques (Belabbas & Chen, 2023).

Notably, the observers of the Koopmanized controlled system received some attentions in literature. Despite the longstanding historic account of mathematical theory of the Koopman operator as well as existing connecting threads between the Koopman operator and dynamic systems in the sense of bilinearizations, there is no publication available yet on control and filtering of the nonlinear Itô stochastic differential system in the Koopman operator framework. Unifying the Koopman theory and filtering of the controlled nonlinear Itô stochastic differential system becomes quite hard. This seems to be one of the reasons that the filtering problem concerning the controlled nonlinear Itô differential system in the Koopman context is yet to be resolved. However, this paper achieves that.

Concerning the Koopmanization and filtering, we sketch the proofs of two Theorems. The proofs are formal, systematic, and rigorous via unifying the notions of functional analysis (Kolmogorov & Fomin, 1957; Limaye, 2016) in conjunction with the Koopman operator theory (Koopman, 1931, Mauroy, Mezić & Sususki, 2020; Brunton, Budišić, Kaiser & Kutz, 2022), Itô calculus (Karatzas & Shreve, 1987) as well as generalized linear system theory. The action of the infinitesimal stochastic Koopman operator on the observable state space encompasses the action on the observable corresponding to the drift, diffusion coefficient and random forcing term as a consequence of the Itô differential rule (Črnjarić-žic, Maćešić & Mezić, 2020).

The organization of the paper is as follows: Section 2 encompasses the main results, i.e., Koopmanization of the controlled nonlinear Itô SDE and its filtering. Further, the main results are rephrased in the Carleman sense as well. Section 3 explains procedure for bilinearizing the controlled nonlinear Itô SDE in the Koopman framework and achieving its filtering. Numerical simulations for the comparison between both the methods can be found in Section 4, whereas concluding remarks are given in Section 5.

## 2. Main results

This section sketches a proof of the 'Koopmanization' concerning the original controlled nonlinear stochastic system in a formal, systematic, and rigorous setting. It is followed by another Theorem and its proof on the Koopman filter.

**Theorem 1.** Consider a finite dimensional controlled nonlinear Itô stochastic differential system $\sum_{1}^{s}$,

$$dx_t = f(x_t, u_t)dt + \sum_{1 \leq \gamma \leq r} g_\gamma(x_t, u_t)dB_\gamma, \quad (1)$$

where the state vector $x_t \in R^n$, the drift vector field $f: R^n \times U \to R^n$, the input control parameter vector $u_t \in U \subset R^d$ and the process noise coefficient $g_\gamma: R^n \times U \to R^n$ with $1 \leq \gamma \leq r$. Given probability space $\{\Omega, \Im, \mu\}$, the stochastic process $X = \{x_t, \Im_t, 0 \le t < \infty\}$ so is $\{(t, \omega) | x_t(\omega) \in R^n, R^n \in B(R^n)\} \subset [0, T] \times \Omega$. The process $\{B_t, \Im_t, 0 \le t < \infty\}$ is a $r$-dimensional standard Brownian motion process on the probability space $\{\Omega, \Im, \mu\}$.

Then, the finite dimensional Koopmanization of the given original system $\sum_1^s$ boils down to the system $\sum_2^s$

$$dz_t = \Lambda_{f'}(u_t) z_t dt + \sum_{1 \le \gamma \le r} D_\gamma(u_t) z_t dB_\gamma + \sum_{1 \le \gamma \le r} F_\gamma(u_t) dB_\gamma \qquad (2)$$

$$= (\Lambda_f(u_t) + \Lambda_b(u_t)) z_t dt + \sum_{1 \le \gamma \le r} D_\gamma(u_t) z_t dB_\gamma + \sum_{1 \le \gamma \le r} F_\gamma(u_t) dB_\gamma,$$

where $\Lambda_{f'}$ is the state matrix associated with the modified drift vector field $f'$ such that $\Lambda_{f'}(u_t) = \Lambda_f(u_t) + \Lambda_b(u_t)$. That accounts for the state matrices $\Lambda_f(u_t)$ and $\Lambda_b(u_t)$ corresponding to the drift and diffusion coefficients of the Itô stochastic differential system respectively.

Select the finite number $m$ such that an error bound for the Koopman operator framework is

$$E \left\| L_{g_\gamma} \varphi(x_t) - \sum_{j=1}^\infty \varphi_j(x_t) v_j^{g_\gamma}(u_t) - F_\gamma(u_t) \right\|$$

$$< \inf_{m < \infty} E \left\| L_{g_\gamma} \varphi(x_t) - \sum_{1 \le j \le m} \varphi_j(x_t) v_j^{g_\gamma}(u_t) - F_\gamma(u_t) \right\| < \varepsilon,$$

for each $\varepsilon > 0$ with all time $t \ge 0$, $z_t \in R^m$, $D_\gamma(u_t) \in R^{m \times m}$, $F_\gamma(u_t) \in R^m$, where $E$ is the conditional expectation operator.

**Proof.** Suppose $n$ is the dimensionality of the Itô stochastic differential system under considerations, $d$ is the dimensionality of the control parameter and $r$ is the dimensionality of the random process acting on the Itô stochastic differential system. A proof of the Theorem concerning the system can be sketched using the action of the Koopman operator on the Hilbert space. Consider a basis $\{\varphi_1, \varphi_2 \ldots, \varphi_m\}$ is a finite subset of $\{\varphi_1, \varphi_2, \ldots\}$. The basis $\{\varphi_1, \varphi_2, \ldots\}$ spans the infinite dimensional eigenfunction space of the Koopman operator. The infinite dimensional eigenfunction space is a subspace of the Hilbert space. The basis $\{\varphi_1, \varphi_2, \ldots, \varphi_m\}$ is an eigenfunction basis of the Koopman operator. The basis $\{\varphi_1, \varphi_2, \ldots, \varphi_m\}$ spans the finite dimensional eigenfunction space of the Koopman operator related to the drift vector field. Invoke the invariant subspace assumption of the Koopman operator related to the process noise coefficient vector field. Suppose the function $\varphi_j(x_t)$ has the analyticity and possess the nice properties that allow the applicability of the multi-dimensional Itô stochastic differential rule admitting the notion of Lie derivatives. Thus, the Lie derivative of $\varphi_j(x_t)$ with $j \ge 1$ is

$$d\varphi_j(x_t) = \sum_{1 \le l \le n} \frac{\partial \varphi_j(x_t)}{\partial x_l} f_l(x_t, u_t) dt + \sum_{1 \le l \le n, 1 \le \gamma \le r} \frac{\partial \varphi_j(x_t)}{\partial x_l} g_\gamma^{(l)}(x_t, u_t) dB_\gamma$$

$$+ \frac{1}{2} \sum_{1 \le l_1 \le n, 1 \le l_2 \le n} \frac{\partial^2 \varphi_j(x_t)}{\partial x_{l_1} \partial x_{l_2}} dx_{l_1} dx_{l_2}. \qquad (3)$$

Now, concerning the setting of the Itô stochastic differential equation of Theorem 1, we have

$$dx_t dx_t^T = \sum_{1\leq \gamma \leq r} g_\gamma(x_t, u_t) g_\gamma^T(x_t, u_t) dt.$$

As a result of the above, the $(l_1, l_2)^{th}$ component of the matrix term $dx_t dx_t^T$ becomes

$$dx_{l_1} dx_{l_2} = (dx_t dx_t^T)_{l_1 l_2} = \left(\sum_{1\leq \gamma \leq r} g_\gamma(x_t, u_t) g_\gamma^T(x_t, u_t) dt\right)_{l_1 l_2}$$
$$= \sum_{1\leq \gamma \leq r} (g_\gamma(x_t, u_t) g_\gamma^T(x_t, u_t))_{l_1 l_2} dt$$
$$= \sum_{1\leq \gamma \leq r} g_\gamma^{(l_1)}(x_t, u_t) g_\gamma^{(l_2)}(x_t, u_t) dt. \quad (4)$$

Equation (3) in conjunction with equation (4), we have

$$d\varphi_j(x_t) = \sum_{1\leq l\leq n} \frac{\partial \varphi_j(x_t)}{\partial x_l} f_l(x_t, u_t) dt + \sum_{1\leq l\leq n, 1\leq \gamma \leq r} \frac{\partial \varphi_j(x_t)}{\partial x_l} g_\gamma^{(l)}(x_t, u_t) dB_\gamma$$
$$+ \frac{1}{2} \sum_{1\leq l_1 \leq n, 1\leq l_2 \leq n} \sum_{1\leq \gamma \leq r} \frac{\partial^2 \varphi_j(x_t)}{\partial x_{l_1} \partial x_{l_2}} g_\gamma^{(l_1)}(x_t, u_t) g_\gamma^{(l_2)}(x_t, u_t) dt$$
$$= \left(\sum_{1\leq l\leq n} \frac{\partial \varphi_j(x_t)}{\partial x_l} f_l(x_t, u_t) + \frac{1}{2} \sum_{1\leq l_1 \leq n, 1\leq l_2 \leq n} \sum_{1\leq \gamma \leq r} \frac{\partial^2 \varphi_j(x_t)}{\partial x_{l_1} \partial x_{l_2}} g_\gamma^{(l_1)}(x_t, u_t) g_\gamma^{(l_2)}(x_t, u_t)\right) dt$$
$$+ \sum_{1\leq l\leq n, 1\leq \gamma \leq r} \frac{\partial \varphi_j(x_t)}{\partial x_l} g_\gamma^{(l)}(x_t, u_t) dB_\gamma.$$

Alternatively,

$$d\varphi_j(x_t) = (L_u^f \varphi_j(x_t) + L_u^b \varphi_j(x_t)) dt + \sum_{1\leq \gamma \leq r} L_u^{g_\gamma} \varphi_j(x_t) dB_\gamma,$$

where

$$L_u^f \varphi_j(x_t) = \left\langle \frac{\partial \varphi_j(x_t)}{\partial x_t}, f(x_t, u_t) \right\rangle, L_u^b \varphi_j(u_t) = \frac{1}{2} \operatorname{Sp}\left\langle \frac{\partial^2 \varphi_j(x_t)}{\partial x_t \partial x_t^T}, \sum_{1\leq \gamma \leq r} g_\gamma(x_t, u_t) g_\gamma^T(x_t, u_t) \right\rangle,$$

$$L_u^{g_\gamma} \varphi_j(x_t) = \left\langle \frac{\partial \varphi_j(x_t)}{\partial x_t}, g_\gamma(x_t, u_t) \right\rangle,$$

and $L$ represents Lie derivative, the superscripts $f$ and $b$ correspond to the drift and diffusion matrix of the Itô stochastic differential system. Now consider $\varphi(x_t) = (\varphi_1(x_t) \quad \varphi_2(x_t) \ldots)^T$ concerning the basis $\{\varphi_1, \varphi_2, \cdots\}$. Thanks to a notion of the Lie derivative as well as a consequence of a closure of the Lie derivative equipped with the operator-invariant subspace property, we have

$$d\varphi(x_t) = (L_u^f \varphi(x_t) + L_u^b \varphi(x_t)) dt + \sum_{1\leq \gamma \leq r} L_u^{g_\gamma} \varphi(x_t) dB_\gamma,$$

$$L_u^{g_\gamma} \varphi(x_t) - F_\gamma(u_t) = \sum_{j\geq 1} \varphi_j(x_t) v_j^{g_\gamma}(u_t).$$

Thus, after combining the above set of equations, we have

$$d\varphi(x_t) = (L_u^f \varphi(x_t) + L_u^b \varphi(x_t))dt + \sum_{1 \leq \gamma \leq r} L_u^{g_\gamma} \varphi(x_t) dB_\gamma$$

$$= (L_u^f \varphi(x_t) + L_u^b \varphi(x_t))dt + \sum_{1 \leq \gamma \leq r} \sum_{j \geq 1} \varphi_j(x_t) v_j^{g_\gamma}(u_t) dB_\gamma + \sum_{1 \leq \gamma \leq r} F_\gamma(u_t) dB_\gamma$$

$$= \left( \sum_{j \geq 1} \varphi_j(x_t) v_j^f(u_t) + \sum_{j \geq 1} \varphi_j(x_t) v_j^b(u_t) \right) dt + \sum_{1 \leq \gamma \leq r} \sum_{j \geq 1} \varphi_j(x_t) v_j^{g_\gamma}(u_t) dB_\gamma$$

$$+ \sum_{1 \leq \gamma \leq r} F_\gamma(u_t) dB_\gamma. \quad (5)$$

Equation (5) is a consequence of the action of the Koopman operator on the infinite-dimensional linear function space. Note that equation (5) is a Koopmanization of the original controlled nonlinear Itô SDE in the infinite-dimensional Koopman operator framework. For the finite dimensional Koopmanization of the nonlinear system, choose the eigenfunction basis $\{\varphi_1, \varphi_2, \cdots, \varphi_m\}$ with the help of operator-invariant subspace property. Thus,

$$d\varphi(x_t) = \left( \sum_{1 \leq j \leq m} \varphi_j(x_t) v_j^f(u_t) + \sum_{1 \leq j \leq m} \varphi_j(x_t) v_j^b(u_t) \right) dt$$

$$+ \sum_{1 \leq \gamma \leq r} \sum_{1 \leq j \leq m} \varphi_j(x_t) v_j^{g_\gamma}(u_t) dB_\gamma + \sum_{1 \leq \gamma \leq r} F_\gamma(u_t) dB_\gamma$$

$$= \left( \Lambda_f(u_t) \varphi(x_t) + \Lambda_b(u_t) \varphi(x_t) \right) dt + \sum_{1 \leq \gamma \leq r} D_\gamma(u_t) \varphi(x_t) dB_\gamma + \sum_{1 \leq \gamma \leq r} F_\gamma(u_t) dB_\gamma, \quad (6)$$

where the terms $\Lambda_f$, $\Lambda_b$, $D_\gamma$ are the matrices, $F_\gamma(u_t)$ is the column vector. After adopting a convenient notation in the augmented observable state space setting associated with equation (6), the system $\sum_2^s$

$$dz_t = (\Lambda_f(u_t) + \Lambda_b(u_t)) z_t dt + \sum_{1 \leq \gamma \leq r} D_\gamma(u_t) z_t dB_\gamma + \sum_{1 \leq \gamma \leq r} F_\gamma(u_t) dB_\gamma, \quad \text{a.s. } (\mu) \quad (7)$$

is a finite-dimensional Koopmanization of the original controlled nonlinear Itô stochastic differential system $\sum_1^s$. Note that $z_t = \varphi(x_t) = (\varphi_j(x_t))_{1 \leq j \leq m}$, $\varphi: R^n \to R^m$ and $m > n$. Equation (7) describes a bilinearization of the Itô SDE in the finite-dimensional Koopman operator framework with the generalized control $u_t$. If the function $u_t$ is non-random and equation possess a unique solution, the associated controlled state is a Markov process.

For the larger $m$, the approximation error arising from the finite dimensional Koopmanization of the original nonlinear system reduces. For the 'deterministic' case, we bound the approximation error attributed to the finite-dimensional Koopman operator framework with the approximation error tends to zero such that $m \to \infty$. We re-formalize it for the 'stochastic' case by using the conditional expectation, i.e.

$$\inf_{1 \leq m} E \left\| L_{g_\gamma} \varphi(x_t) - \sum_{1 \leq j \leq m} \varphi_j(x_t) v_j^{g_\gamma}(u_t) - F_\gamma(u_t) \right\|$$

$$= E \left\| L_{g_\gamma} \varphi(x_t) - \sum_{j=1}^{\infty} \varphi_j(x_t) v_j^{g_\gamma}(u_t) - F_\gamma(u_t) \right\|.$$

The embedding of the greater number of eigenfunctions into the finite dimensional representation of the infinite dimensional Koopman bilinear system circumvents the curse of dimensionality. Thus,

'embedding the greater number of eigenfunctions with the given state at the time $t$ and a choice of the control parameter' for each $0 < \varepsilon < \infty$ that we bound the approximation error in the sense of conditional expectation of stochastic processes such that

$$E \left\| L_{g_\gamma} \varphi(x_t) - \sum_{j=1}^{\infty} \varphi_j(x_t) v_j^{g_\gamma}(u_t) - F_\gamma(u_t) \right\|$$

$$< \inf_{m < \infty} E \left\| L_{g_\gamma} \varphi(x_t) - \sum_{1 \leq j \leq m} \varphi_j(x_t) v_j^{g_\gamma}(u_t) - F_\gamma(u_t) \right\| < \varepsilon$$

for all time $t \geq 0$.

*QED*

Consider the measurement system $dy_t = h(x_t)dt + r_1 d\eta_t$, where $h(x_t)$ is a linear combination of observables in the finite dimensional state space $R^m$. As a result of this,

$$dy_t = h(x_t)dt + r_1 d\eta_t = (\sum_{1 \leq j \leq m} \varphi_j(x_t) v_j^{(h)}) dt + r_1 d\eta_t = C^{(h)} \varphi(x_t) dt + r_1 d\eta_t$$

$$= C^{(h)} z_t dt + r_1 d\eta_t,$$

(8)

where $y_t \in R^p$, $C^{(h)} \in R^{p \times m}$ and the process $\{\eta_t, \Im_t, 0 \leq t < \infty\}$ is a $p$-dimensional standard Brownian motion process on the probability space $\{\Omega, \Im, \mu\}$.

*Koopman-based nonlinear filter*

Filtering is the stochastic terminology corresponding to the observer of dynamic systems. The filtering in the Koopman sense can be weaved by utilizing the stochastic differential system in the Koopman operator framework. The Koopman operator framework of the nonlinear Itô stochastic differential system culminates into 'the bilinear Itô stochastic differential system in the augmented state space'. Theorem 2 of the paper is about a systematic construction of the Koopman-based nonlinear filter with its proof.

**Theorem 2.** Consider a stochastic Koopman system $\sum_2^s$ with the measurement system

$$dz_t = \left( \Lambda_f(u_t) + \Lambda_b(u_t) \right) z_t dt + \sum_{1 \leq \gamma \leq r} D_\gamma(u_t) z_t dB_\gamma + \sum_{1 \leq \gamma \leq r} F_\gamma(u_t) dB_\gamma,$$

$$dy_t = C^{(h)} z_t dt + r_1 d\eta_t,$$

where $z_t \in R^m$, $C^{(h)} \in R^{p \times m}$. Note that the stochastic Koopman system $\sum_2^s$ has the lifted dynamics associated with the original controlled nonlinear Itô differential system $\sum_1^s$. Suppose all conditions and notations of Theorem 1 hold. The Koopman filter $\widehat{\sum}_2^s$ associated with the system $\sum_2^s$ boils down to the following

$$d\hat{z}_t = (\Lambda_f(u_t) + \Lambda_b(u_t))\hat{z}_t dt + P_{z_t} C^{(h)^T} r_1^{-2} (dy_t - C^{(h)} \hat{z}_t dt),$$

$$(dP_{z_t})_{ij} = (P_{z_t}\Lambda_f^T(u_t) + P_{z_t}\Lambda_b^T(u_t) + \Lambda_f(u_t)P_{z_t} + \Lambda_b(u_t)P_{z_t} + \sum_{1\leq\gamma\leq r} D_\gamma(u_t) P_{z_t} D_\gamma^T(u_t)$$

$$+ \sum_{1\leq\gamma\leq r} D_\gamma(u_t)\hat{z}_t\hat{z}_t^T D_\gamma^T(u_t) + \sum_{1\leq\gamma\leq r} D_\gamma(u_t)\hat{z}_t F_\gamma^T(u_t)$$

$$+ \sum_{1\leq\gamma\leq r} F_\gamma(u_t)\hat{z}_t^T D_\gamma^T(u_t) + \sum_{1\leq\gamma\leq r} F_\gamma(u_t)F_\gamma^T(u_t) - P_{z_t}C^{(h)T}r_1^{-2}C^{(h)}P_{z_t})dt$$

such that $\hat{z}_t = E(z_t|y_\tau, 0\leq\tau\leq t)$, $P_{z_t} = E((z_t-\hat{z}_t)(z_t-\hat{z}_t)^T|y_\tau, 0\leq\tau\leq t)$.

**Proof.** Recall the system with control case. The controlled system is identically equal to the system with no control. Thus,

$$dz_t = (\Lambda_f(u_t) + \Lambda_b(u_t))z_t dt + \sum_{1\leq\gamma\leq r} D_\gamma(u_t) z_t dB_\gamma + \sum_{1\leq\gamma\leq r} F_\gamma(u_t)dB_\gamma,$$

$$= (\Lambda_f^{cl} + \Lambda_b^{cl})z_t dt + \sum_{1\leq\gamma\leq r} D_\gamma^{cl} z_t dB_\gamma + \sum_{1\leq\gamma\leq r} F_\gamma^{cl} dB_\gamma, \tag{9a}$$

where $\Lambda_f^{cl} \in R^{m\times m}$, $\Lambda_b^{cl} \in R^{m\times m}$, $D_\gamma^{cl} \in R^{m\times m}$ and $F_\gamma^{cl} \in R^m$. We wish to construct a Koopman based nonlinear filter concerning 'equation (9a) in combination with the linear observation equation $dy_t = C^{(h)}z_t dt + r_1 d\eta_t$.' Since the Koopman-based filtering is a nonlinear filtering, utilize the coupled nonlinear filtering equations (Pugachev & Sinitsyn, 1987) with appropriate notations for the drift, process noise coefficient and measurement nonlinearity, i.e. $a(z_t)$, $b(z_t)$ and $h(z_t)$ respectively, such that

$$a(z_t) = (\Lambda_f^{cl} + \Lambda_b^{cl})z_t, \quad b(z_t)dB_t = \sum_{1\leq\gamma\leq r} D_\gamma^{cl} z_t dB_\gamma + \sum_{1\leq\gamma\leq r} F_\gamma^{cl} dB_\gamma, \quad h(z_t) = C^{(h)}z_t. \tag{9b}$$

Note that $h(z_t) = C^{(h)}z_t$ from the considerations of the linear measurement system. As a result of this, we arrive at the Koopman-based nonlinear filtering equations. An analytically tractable form of the nonlinear filtering (Pugachev & Sinitsyn, 1987) can be found in (Sharma, 2009). A system of nonlinear filtering equations with 'appropriate notations' are

$$d\hat{z}_t = (a(\hat{z}_t) + \frac{1}{2}\sum_{p,q} P_{pq}\frac{\partial^2 a(\hat{z}_t)}{\partial\hat{z}_p\partial\hat{z}_q})dt + P_{z_t}\frac{\partial h^T(\hat{z}_t)}{\partial\hat{z}_t}r_1^{-2}(dy_t - h(\hat{z}_t)dt$$

$$-\frac{1}{2}\sum_{p,q} P_{pq}\frac{\partial^2 h(\hat{z}_t)}{\partial\hat{z}_p\partial\hat{z}_q}dt), \tag{10a}$$

$$(dP_z)_{ij} = (\sum_p P_{ip}\frac{\partial a_j(\hat{z}_t)}{\partial\hat{z}_p} + \sum_p P_{jp}\frac{\partial a_i(\hat{z}_t)}{\partial\hat{z}_p} + (bb^T)_{ij}(\hat{z}_t) + \frac{1}{2}\sum_{p,q} P_{pq}\frac{\partial^2(bb^T)_{ij}(\hat{z}_t)}{\partial\hat{z}_p\partial\hat{z}_q}$$

$$-(\sum_p P_{ip}\frac{\partial h^T(\hat{z}_t)}{\partial\hat{z}_p})r_1^{-2}(\sum_p P_{jp}\frac{\partial h(\hat{z}_t)}{\partial\hat{z}_p}))dt$$

$$+(\sum_{p,q} P_{ip}P_{jq}\frac{\partial^2 h^T(\hat{z}_t)}{\partial\hat{z}_p\partial\hat{z}_q})r_1^{-2}(dy_t - h(\hat{z}_t)dt - \frac{1}{2}\sum_{p,q} P_{pq}\frac{\partial^2 h(\hat{z}_t)}{\partial\hat{z}_p\partial\hat{z}_q}dt). \tag{10b}$$

The set of equations (10a) - (10b) exploits the nonlinear Itô SDE setting with the state $z_t \in R^m$ of equation (9) in conjunction with equation (8). Note that the state $z_t$ is an augmented state vector associated with the Koopmanized SDE (7) related to equation (1). Using equation (9b), we get a useful result for the Koopmanization, i.e.

$$b(z_t)b^T(z_t) = \sum_{1 \leq \gamma \leq r} D_\gamma^{cl} z_t z_t^T D_\gamma^{cl\,T} + \sum_{1 \leq \gamma \leq r} D_\gamma^{cl} z_t F_\gamma^{cl\,T} + \sum_{1 \leq \gamma \leq r} F_\gamma^{cl} z_t^T D_\gamma^{cl\,T} + \sum_{1 \leq \gamma \leq r} F_\gamma^{cl} F_\gamma^{cl\,T}$$

where $(bb^T)_{ij}(z_t, t) = (b(z_t, t)b^T(z_t, t))_{ij}$. After the action of the conditional expectation, given noisy observables accumulated upto the time $t$ to both the sides of above equation, we have

$$E(b(z_t)b^T(z_t)|Y_t) = \sum_{1 \leq \gamma \leq r} D_\gamma^{cl} P_{z_t} D_\gamma^{cl\,T} + \sum_{1 \leq \gamma \leq r} D_\gamma^{cl} \hat{z}_t \hat{z}_t^T D_\gamma^{cl\,T} + \sum_{1 \leq \gamma \leq r} D_\gamma^{cl} \hat{z}_t F_\gamma^{cl\,T}$$

$$+ \sum_{1 \leq \gamma \leq r} F_\gamma^{cl} \hat{z}_t^T D_\gamma^{cl\,T} + \sum_{1 \leq \gamma \leq r} F_\gamma^{cl} F_\gamma^{cl\,T}, \qquad (11)$$

where $E(z_t z_t^T | y_\tau, 0 \leq \tau \leq t) = P_{z_t} + \hat{z}_t \hat{z}_t^T$. Equation (11) in conjunction with equations (10a)-(10b) and (9b), we have the following Koopman-based filtering equations.

$$d\hat{z}_t = (\Lambda_f^{cl} + \Lambda_b^{cl})\hat{z}_t dt + P_{z_t} C^{(h)T} r_1^{-2}(dy_t - C^{(h)} \hat{z}_t dt), \qquad (12a)$$

$$dP_{z_t} = (P_{z_t} \Lambda_f^{cl\,T} + P_{z_t} \Lambda_b^{cl\,T} + \Lambda_f^{cl} P_{z_t} + \Lambda_b^{cl} P_{z_t} + \sum_{1 \leq \gamma \leq r} D_\gamma^{cl} P_{z_t} D_\gamma^{cl\,T} + \sum_{1 \leq \gamma \leq r} D_\gamma^{cl} \hat{z}_t \hat{z}_t^T D_\gamma^{cl\,T}$$

$$+ \sum_{1 \leq \gamma \leq r} D_\gamma^{cl} \hat{z}_t F_\gamma^{cl\,T} + \sum_{1 \leq \gamma \leq r} F_\gamma^{cl} \hat{z}_t^T D_\gamma^{cl\,T} + \sum_{1 \leq \gamma \leq r} F_\gamma^{cl} F_\gamma^{cl\,T}$$

$$- P_{z_t} C^{(h)T} r_1^{-2} C^{(h)} P_{z_t}) dt. \qquad (12b)$$

It is worth to mention that equation (11) makes the conditional variance evolution expressible in the matrix-vector format. The 'generalized' Riccati equation combines 'the Riccati equation and the corrections attributed to the diffusion coefficient'. The generalized Riccati equation with 'the matrix $D_\gamma^{cl}$ is identically zero with every $\gamma'$ reduces to the Riccati equation. Thus, the conditional variance evolution of the Koopmanized SDE is a generalized Riccati equation, i.e., equation (12b). Recall equation (9a), then the "controlled version of 'the nonlinear filtering (12a)-(12b) in the Koopman sense' is the system of the following coupled equations":

$$d\hat{z}_t = (\Lambda_f(u_t) + \Lambda_b(u_t))\hat{z}_t dt + P_{z_t} C^{(h)T} r_1^{-2}(dy_t - C^{(h)} \hat{z}_t), \qquad (13a)$$

$$(dP_{z_t})_{ij} = (P_{z_t} \Lambda_f^T(u_t) + P_{z_t} \Lambda_b^T(u_t) + \Lambda_f(u_t) P_{z_t} + \Lambda_b(u_t) P_{z_t} + \sum_{1 \leq \gamma \leq r} D_\gamma(u_t) P_{z_t} D_\gamma^T(u_t)$$

$$+ \sum_{1 \leq \gamma \leq r} D_\gamma(u_t) \hat{z}_t \hat{z}_t^T D_\gamma^T(u_t) + \sum_{1 \leq \gamma \leq r} D_\gamma(u_t) \hat{z}_t F_\gamma^T(u_t)$$

$$+ \sum_{1 \leq \gamma \leq r} F_\gamma(u_t) \hat{z}_t^T D_\gamma^T(u_t) + \sum_{1 \leq \gamma \leq r} F_\gamma(u_t) F_\gamma^T(u_t) - P_{z_t} C^{(h)T} r_1^{-2} C^{(h)} P_{z_t}) dt. \qquad (13b)$$

**Remark 1.** Consider the finite dimensional Carleman-embedded controlled stochastic differential system

$$d\xi_t = \Lambda(u_t) \xi_t dt + \sum_{1 \leq \gamma \leq r} D_\gamma(u_t) \xi_t dB_\gamma + \sum_{1 \leq \gamma \leq r} F_\gamma(u_t) dB_\gamma, \qquad \text{a.s. } (\mu)$$

with the measurement system $dy_t = C\xi_t dt + r_1 d\eta_t$. Similar to the construction of the proof of Theorem 2, a system of nonlinear filtering equations in the Carleman sense concerning 'the controlled lifted stochastic dynamics of equation (6) in combination with linear measurements' boils down to

$$d\hat{\xi}_t = \Lambda(u_t)\hat{\xi}_t dt + P_{\xi_t} C^T r_1^{-2}(dy_t - C\hat{\xi}_t dt), \tag{14a}$$

$$dP_{\xi_t} = (P_{\xi_t}\Lambda^T(u_t) + \Lambda(u_t)P_{\xi_t}(u_t) + \sum_{1\leq\gamma\leq r} D_\gamma(u_t) P_{\xi_t} D_\gamma^T(u_t)$$

$$+ \sum_{1\leq\gamma\leq r} D_\gamma(u_t) \hat{\xi}_t\hat{\xi}_t^T D_\gamma^T(u_t) + \sum_{1\leq\gamma\leq r} D_\gamma(u_t)\hat{\xi}_t F_\gamma^T(u_t)$$

$$+ \sum_{1\leq\gamma\leq r} F_\gamma(u_t)F_\gamma^T(u_t) - P_{\xi_t} C^T r_1^{-2} C P_{\xi_t})dt. \tag{14b}$$

Note that the notations of the matrices $\Lambda(u_t), D_\gamma(u_t), F_\gamma(u_t)$ and $P_{\xi_t}$ associated with the filtering in the Carleman sense, i.e. (14a)-(14b), have different interpretations and different sizes from that of the filtering in the Koopman sense given in equations (13a)-(13b). The Koopmanization of the original system hinges on the operator-theoretic framework; on the other hand, the Carleman linearization exploits the linearization about a given point of the system trajectory. Importantly, both the methods adopt different frameworks, but they do the bilinearization of the controlled nonlinear system. Thus, equation (14b) is a generalized Riccati equation in the Carleman sense, which is different from the generalized Riccati equation in the Koopman sense, see equation (13b).

## 3. Illustrations of the main results: Koopman framework

This section is concerning the bilinearization and filtering of the nonlinear Itô SDE in Koopman framework. Consider a polynomial nonlinear system in the Itô setting (Germani, Manes & Palumbo, 2007):

$$dx_1 = (-x_1 + x_1 x_2)dt + adB_t, \quad dx_2 = (-2x_2 - 2x_1 x_2)dt + bdB_t, \tag{15a}$$

$$dy = (x_1 - x_1 x_2)dt + r_1 d\eta_t, \tag{15b}$$

where $B_t$ and $\eta_t$ are the Brownian motion processes.

In order to bilinearize the above nonlinear system, we first need to find the principal eigenfunction for which consider only autonomous part of the system dynamics. The principal eigenfunction for given system is as follows:

$$\varphi = 2x_1 + x_2 + 2\log x_1 - \log x_2. \tag{16}$$

The generalized eigenfunction state vector for the given nonlinear system is as follows:

$$z_t = (z_1 \quad z_2 \quad z_3 \quad z_4 \quad z_5 \quad z_6)^T = (2x_1 + x_2 + 2\log x_1 - \log x_2 \quad x_1 \quad x_2 \quad x_1^2 \quad x_2^2 \quad x_1 x_2)^T$$

and their system of Itô stochastic evolutions are

$$dz_1 = (\frac{8a^2}{x_{01}^3}z_2 - \frac{4b^2}{x_{02}^3}z_3 - \frac{3a^2}{x_{01}^4}z_4 + \frac{3b^2}{2x_{02}^4}z_5)dt + (2a + b + \frac{6a}{x_{01}} - \frac{3b}{x_{02}}$$

$$- \frac{6a}{x_{01}^2}z_2 + \frac{3b}{x_{02}^2}z_3 + \frac{2a}{x_{01}^3}z_4 - \frac{b}{x_{02}^3}z_5)dB_t,$$

$$dz_2 = -(z_2 + z_6)dt + adB_t, dz_3 = (-2z_3 - 2z_6)dt + bdB_t,$$
$$dz_4 = (a^2 - 2z_4)dt + 2az_2dB_t, dz_5 = (b^2 - 4z_5)dt + 2bz_3dB_t,$$
$$dz_6 = -3z_6 dt + (az_3 + bz_2)dB_t. \tag{17a}$$

The output equation becomes
$$dy = (z_2 - z_6)dt + r_1 d\eta_t. \tag{17b}$$

Alternatively, equation (17a) corresponding to the system of equations (15a) can be recast as
$$dz_t = \Lambda_{f'} z_t dt + Dz_t dB_t + F dB_t. \tag{18a}$$

Analogously, equation (17b) corresponding to equation (15b) can be recast as
$$dy_t = C^{(h)} z_t dt + r_1 d\eta_t, \tag{18b}$$

where

$$\Lambda_{f'} = \begin{pmatrix} 0 & \frac{8a^2}{x_{01}^3} & -\frac{4b^2}{x_{02}^3} & -\frac{3a^2}{x_{01}^4} & \frac{3b^2}{2x_{02}^4} & 0 \\ 0 & -1 & 0 & 0 & 0 & 1 \\ 0 & 0 & -2 & 0 & 0 & -2 \\ 0 & 0 & 0 & -2 & 0 & 0 \\ 0 & 0 & 0 & 0 & -4 & 0 \\ 0 & 0 & 0 & 0 & 0 & -3 \end{pmatrix}, F = \begin{pmatrix} 2a + b + \frac{6a}{x_{01}} - \frac{3b}{x_{02}} \\ a \\ b \\ 0 \\ 0 \\ 0 \end{pmatrix},$$

$$D = \begin{pmatrix} 0 & -\frac{6a}{x_{01}^2} & \frac{3b}{x_{02}^2} & \frac{2a}{x_{01}^3} & -\frac{b}{x_{02}^3} & 0 \\ 0 & 0 & 0 & 0 & 0 & 0 \\ 0 & 0 & 0 & 0 & 0 & 0 \\ 0 & 2a & 0 & 0 & 0 & 0 \\ 0 & 0 & 2b & 0 & 0 & 0 \\ 0 & b & a & 0 & 0 & 0 \end{pmatrix}, C^{(h)} = (0 \ 1 \ 0 \ 0 \ 0 \ -1).$$

The subscript notation $f'$ recalls the combined contribution of the drift and diffusion coefficient of the Itô SDE (15a) that is attributed to the multi-dimensional Itô differential rule.

*Koopman Filtering*

Filtering of 'the SDE with observation equation, i.e. (18a) and (18b), in the Koopman setting' is the direct consequence of equation (13a) of the proof of Theorem 2 of the paper. Thus,

$$d\hat{z}_1 = \left( \frac{8a^2}{x_{01}^3} \hat{z}_2 - \frac{4b^2}{x_{02}^3} \hat{z}_3 - \frac{3a^2}{x_{01}^4} \hat{z}_4 + \frac{3b^2}{2x_{02}^4} \hat{z}_5 \right) dt$$
$$+ (P_{z_1 z_2} - P_{z_1 z_6}) r_1^{-2} (dy_t - (\hat{z}_2 - \hat{z}_6)dt), \tag{19a}$$

$$d\hat{z}_2 = (-\hat{z}_2 + \hat{z}_6)dt + (P_{z_2 z_2} - P_{z_2 z_6}) r_1^{-2} (dy_t - (\hat{z}_2 - \hat{z}_6)dt), \tag{19b}$$

$$d\hat{z}_3 = (-2\hat{z}_3 - 2\hat{z}_6)dt + (P_{z_3 z_2} - P_{z_3 z_6}) r_1^{-2} (dy_t - (\hat{z}_2 - \hat{z}_6)dt), \tag{19c}$$

$$d\hat{z}_4 = -2\hat{z}_4 dt + (P_{z_4 z_2} - P_{z_4 z_6}) r_1^{-2} (dy_t - (\hat{z}_2 - \hat{z}_6)dt), \tag{19d}$$

$$d\hat{z}_5 = -4\hat{z}_5 dt + (P_{z_5 z_2} - P_{z_5 z_6}) r_1^{-2} (dy_t - (\hat{z}_2 - \hat{z}_6)dt), \tag{19e}$$

$$d\hat{z}_6 = -3\hat{z}_6 dt + (P_{z_6 z_2} - P_{z_6 z_6}) r_1^{-2} (dy_t - (\hat{z}_2 - \hat{z}_6)dt). \tag{19f}$$

The concerning generalized Riccati equation in the Koopman sense, i.e. the conditional variance matrix differential equation, can be obtained by using 'equation (13b) of the proof of Theorem 2 of the paper'.

**4. Numerical analysis: Koopman vs. Carleman in stochastic case**

This section demonstrates the numerical simulations of the Koopmanized stochastic differential equation as well as Carleman linearized stochastic differential equation concerning the given nonlinear polynomial system given in equation (15a). First, the efficacy of Koopman bilinearization technique is analysed and then the performance of the corresponding filtering is compared with the filtering in the Carleman framework.

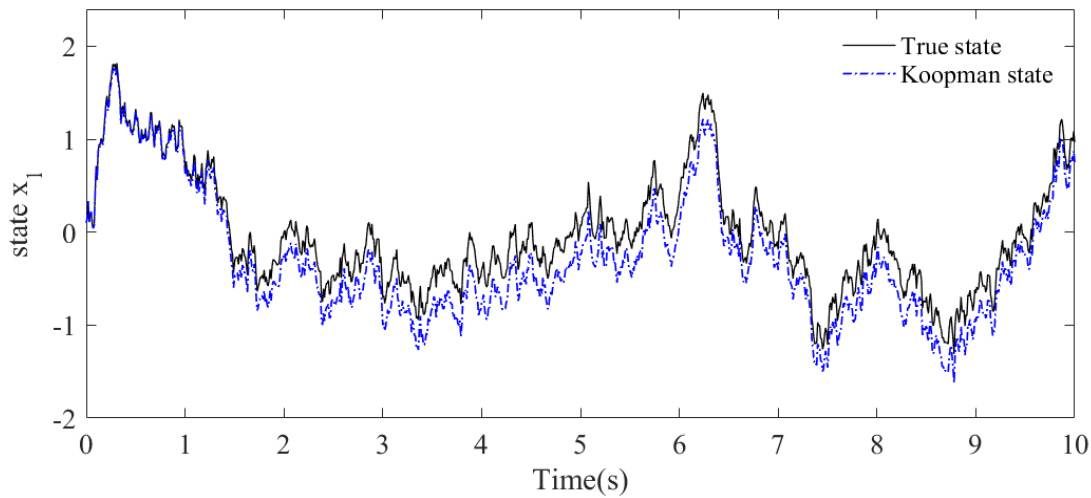

Fig.1. Koopman bilinearized state $x_1$

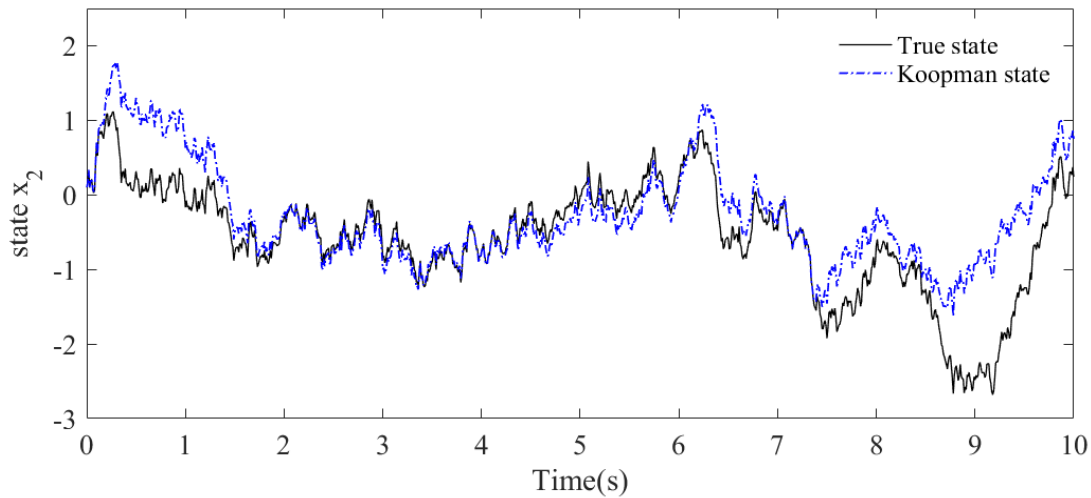

Fig.2. Koopman bilinearized state $x_2$

The numerical simulations are executed using MATLAB© software on Intel(R) Core (TM) i5-8265U laptop, CPU clocked at 1.60GHz with 8.00GB RAM. Choose the system parameters $a = 0.5$, $b = 0.5$, $r_1 = 0.5$ such that the trajectory of the given system is controlled. Choose the initial conditions $x_1(0) = 0.1$, $x_2(0) = 0.1$ and the linearization points of the principal eigenfunction are $x_{01} = -1$ and $x_{02} = 1$ that are admissible with the logarithmic function associated with the principal eigenfunction of equation (16) related to the given polynomial system.

Fig. 1 and Fig. 2 show comparison between the two trajectories, the Koopman bilinearized trajectory and corresponding actual state trajectory. The former is associated with the state $x_1$ of the original polynomial system and the latter is associated with the state $x_2$. The randomness in the Koopmanized state trajectory is attributed to the stochastic part associated with equation (18a), on the other, the stochastic term of equation (15a) contributes to the randomness in the actual state trajectory. Interestingly, the Koopman bilinearized trajectories are sufficiently well-behaved in the sense of becoming agreeable with the increase and decrease in the corresponding actual trajectories. That unfold the usefulness of the 'Koopmanization' of the original nonlinear polynomial stochastic differential system from the filtering and control perspectives.

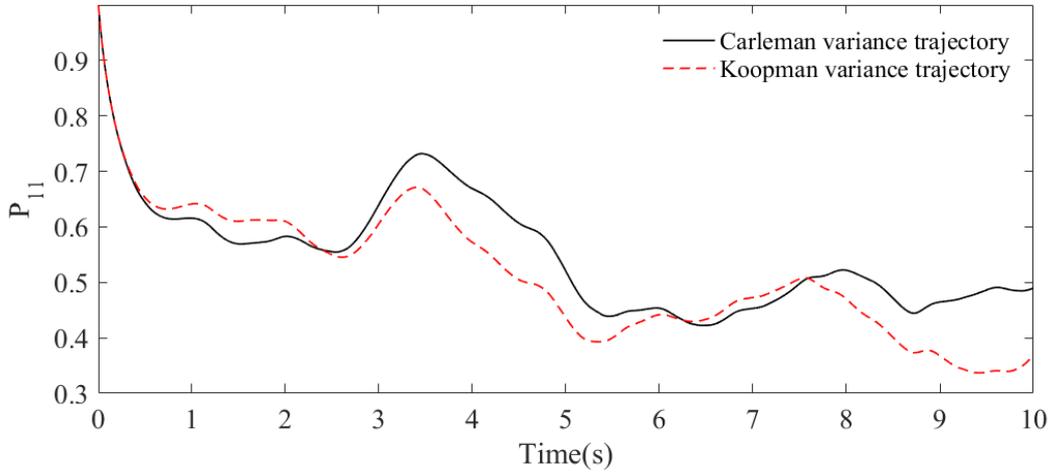

Fig. 3. Conditional variance trajectory $P_{11}$

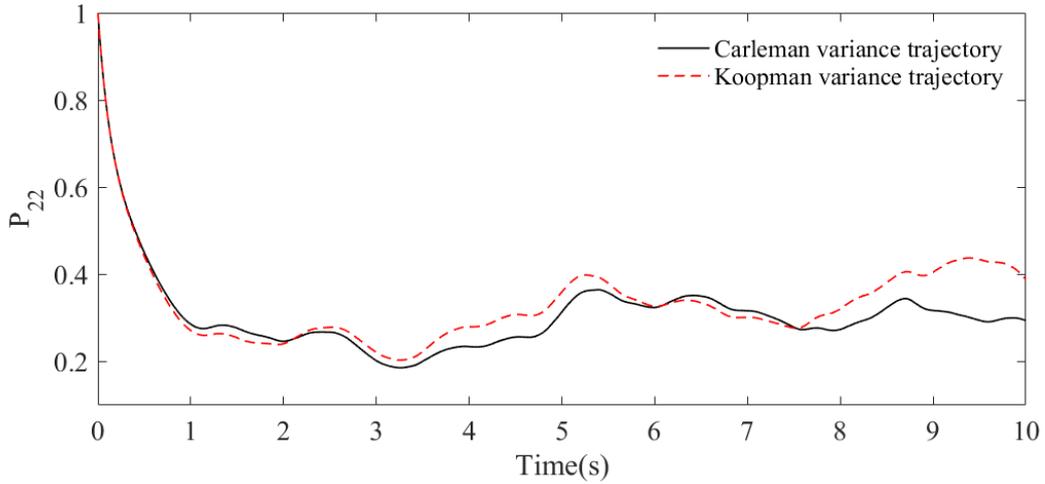

Fig. 4. Conditional variance trajectory $P_{22}$

Further, the comparative analysis is performed between the Koopmans filtering technique and the Carleman filtering technique. The filtering simulation parameters are $\hat{x}_1(0) = 0.1$, $\hat{x}_2(0) = 0.1$, $P_z(0) = I_{6 \times 6}$ and $P_\xi(0) = I_{5 \times 5}$. The Koopmanized filtered state trajectories obey the set of equation (19). The nonlinear filtering accounts for contributions stemming from the system dynamics part and the noisy observations part. For the specific problem considered here, we explain a property of the filtered state trajectory concerning the observation part. For the 'smaller' observation noise intensity, the filtered state estimates in the Koopman sense exhibit small increase or decrease in their trajectory, on the other hand, the filtered estimates in the Carleman framework exhibit relatively larger increase or decrease. Thus, the Koopman framework produces nicer filtered state trajectories than the Carleman. Fig. 3 and Fig. 4 show the conditional variance trajectories for both the frameworks. We infer the less random fluctuations associated with the Koopman variance trajectory in contrast to the Carleman.

The Koopmanization is an operator-theoretic framework for the bilinearization without explicit linearizations. On the other hand, the 'Carleman linearization' of the controlled nonlinear system is a bilinearization, which is expressible in the Kronecker product with the lexicographic ordering arising from the linearization of the nonlinear system combined with the nonlinearity state augmentation. Thus, the Koopman framework has relatively less computational time in contrast to the Carleman linearization. For the given polynomial system of this paper, the computational time for a single simulation run in the Koopman setting is 0.069 secs and 0.077 secs with the Carleman.

## 5. Conclusion

The main theoretical contributions of the paper are the construction of formal, unified, and rigorous proofs of two Theorems concerning the Koopmanization of the controlled Markov processes nonlinear Itô stochastic differential system and then, its filtering. The proofs of the concerning Theorems hinge on the unification of the notions of functional analysis with the linear Koopman operator framework, Itô theory of stochastic processes and the generalized linear system-theoretic framework.

This paper is the first paper of its kind that achieves the Koopmanization via bilinearizing the controlled nonlinear Itô stochastic differential system and then, accomplishes its filtering. This paper explains explicitly "how the 'Koopmanization' is a potential alternative to the 'Carleman linearization' in the sense of bilinearizing and filtering in the Itô stochastic contexts".